# Risk-Averse Joint Capacity Evaluation of PV Generation and Electric Vehicle Charging Stations in Distribution Networks

Huimiao Chen and Zechun Hu
Department of Electrical Engineering
Tsinghua University
Beijing, China
chenhm15@mails.tsinghua.edu.cn

Yinghao Jia
Department of Industrial Engineering
Tsinghua University
Beijing, China

Zuo-Jun Max Shen
Department of Industrial Engineering and Operations Research and Department of Civil and Environmental Engineering
University of California, Berkeley
Berkeley, CA, USA

*Abstract*—**Increasing penetration of distribution generation (DG) and electric vehicles (EVs) calls for an effective way to estimate the achievable capacity connected to the distribution systems, but the exogenous uncertainties of DG outputs and EV charging loads make it challengeable. This study provides a joint capacity evaluation method with a risk threshold setting function for photovoltaic (PV) generation and EV charging stations (EVCSs). The method is mathematically formulated as a distributionally robust joint chance constrained programming model. And the worst-case conditional value at risk (WC-CVaR) approximation and an iterative algorithm based on semidefinite program (SDP) are used to solve the model. Finally, the method test is carried out numerically on IEEE 33-bus radial distribution system.**

## I. INTRODUCTION

In U.S., 30% greenhouse gas emissions come from coal-fired and gas-fired electricity production [1], [2]. As an immediate method towards reducing dependency on fossil fuels and easing environmental pressure, the supply of renewable energy has aroused worldwide attention [3]. In consequence, distributed generation (DG), such as photovoltaic (PV) generation, has already received a great deal of support from political policies and been regarded as a significant portion in the future power system. In addition to the advantage of sustainability, compared with traditional centralized power plants, DG plants are located directly in distribution systems and can avoid the network losses due to the transmission of electricity from a generator to a typical end-user, which accounts for 4.2% to 8.9% of the total electric energy [4]. For almost the same purpose of mitigating energy and environmental issues, electric vehicles (EVs) are expected to take place of internal combustion engine vehicles and accelerate the transformation of the transportation energy of which the global warming gas emission proportion is only second to that of electricity production [5]. Governments around the world have released incentive policies, such as tax credits, rebates and grants, to launch EV demonstration projects, promote the adoption of EVs and cultivate new auto markets [6], [7].

In spite of the above motivations and actions, the advent of DG and EVs calls for new planning and designing techniques in distribution networks. Barriers to integrate DG and EVs into distribution networks include: 1) uncertain DG outputs may cause power security, reliability and quality problems, such as reverse power flow and overvoltage [8]; 2) charging loads of large-scale EV fleets may increase network losses [9], overloading of transformers [10], and excessively heavy line loads [11], etc. Thus, capacity evaluation and optimization of DG generation and EV charging stations (EVCSs) in distribution networks are crucial for the planning.

In existing literature on the above issues, the research of DG and EVCSs are mostly done separately. Strategies and analysis, focusing on the insertion of DG into distribution networks, provide solutions in many respects. For the sake of minimizing the network losses, reference [12] investigates the planning problem of multiple DG units; for the maximum DG penetration, reference [13] tackles the DG integration problem considering harmonic distortion and protection coordination; multi-objective optimization, including DG costs, active and reactive power losses and voltage profiles, of DG placement are carried out in [14]-[16]. Concerning the planning and assessment methods for EVCSs, reference [17] formulates a robust optimization planning methodological framework with the constraints of the power system and the transportation for sustainable integration of plug-in hybrid EVs into a power system; Wang *et al.* [18] develop a multi-objective EVCS planning method which can ensure charging service while reducing energy losses and voltage fluctuations of distribution networks; reference [9] provides a comprehensive approach for evaluating the impact of different levels of EV penetration on distribution network investment and incremental energy losses. As for the research considering both DG and EVCSs, authors

This work was supported in part by the National Natural Science Foundation of China (51477082) and the National Key Research and Development Program (2016YFB0900103).

of [19] propose a multi-year multi-objective planning algorithm for enabling distribution networks to accommodate high penetrations of EVs in conjunction with DG; in [20], a method is developed to obtain the optimal siting and sizing of EVCSs and DG; reference [21] presents an analytical approach to determine the size of EVCSs powered by grid-connected PV penetration; a two-stage approach for allocation of EVCSs and DG in distribution networks is proposed in [22] considering both the economic benefits and network constraints.

Despite the already done work, the forecasting of PV outputs and EV charging load demands is still challengeable, and the uncertainties in the presence of both PV generation and EVCSs have not been investigated adequately. Additionally, as two currently predominant methods to deal with uncertainty, the robust optimization focuses on the worst case in the predetermined ranges of random variables and thereby gives rise to over-conservative results [17], and the stochastic optimization usually requires specific probability distribution functions of random variables which may be unavailable [23]. Hence, a joint capacity evaluation method for PV generation and EVCSs with less information of uncertainties and more careful consideration of risk is necessary for practical planning in distribution networks. To this end, in this paper, a risk-averse joint capacity evaluation method is proposed to deal with the problem. The main procedures and contributions of the paper are summarized below.

- Formulate a distributionally robust (DR) joint chance constrained programming model for joint capacity evaluation of PV generation and EVCSs, which takes full advantage of the easily accessible data, i.e., moment and support sets, of uncertain parameters to describe the uncertainties with a risk budget setting.
- Utilize the worst-case conditional value at risk (WC-CVaR) approximation to reformulate the intractable DR joint chance constraints as a DR individual chance constraint.
- Transform the DR individual chance constraint into approximate bilinear matrix inequality (BMI) constraints.
- Construct an auxiliary optimization problem and develop an iterative algorithm based on semidefinite program (SDP) to solve the BMI constrained problem.

The remainder of the paper is organized as follows. Section II illustrates the mathematical formulation of the evaluation method. Section III provides the solution to the model and Section IV shows the numerical tests. Section V concludes.

## II. MATHEMATICAL FORMULATION OF RISK-AVERSE JOINT CAPACITY EVALUATION

### A. Formulation of Model without Uncertainty

As the basis for the subsequent DR chance constrained programming model, herein, we first introduce a foundational capacity evaluation model without uncertainty. The model aims to obtain the maximum achievable joint capacity of PV generation and EVCSs in the distribution networks, and the objective function is written in weighted form as below.

$$\min \sum_{i \in \Phi^{PV}} \alpha_i^{PV} S_i^{PV} + \sum_{i \in \Phi^{EV}} \alpha_i^{EV} S_i^{EV} \quad (1)$$

where $S_i^{PV}$ and $S_i^{EV}$ are the achievable capacities of PV generation and EVCSs at bus $i$, respectively, and $\alpha_i^{PV}$ and $\alpha_i^{EV}$ are the corresponding weighting coefficients. $\Phi^{PV}$ and $\Phi^{EV}$ are the sets of candidate buses for PV generators and EVCSs, respectively. The constraints are presented as follows.

$$P_{i,k}^{PV} - P_{i,k}^{EV} - P_{i,k}^{OT} - \sum_{j \in \Omega_i} p_{ij,k} = 0, \forall i \in \Phi, k \in \mathbf{T} \quad (2)$$

$$Q_{i,k}^{PV} - Q_{i,k}^{EV} - Q_{i,k}^{OT} - \sum_{j \in \Omega_i} q_{ij,k} = 0, \forall i \in \Phi, k \in \mathbf{T} \quad (3)$$

$$V_{i,k}^2 - V_{j,k}^2 - 2 \cdot (r_{ij} \cdot p_{ij,k} + x_{ij} \cdot q_{ij,k}) - \frac{(r_{ij}^2 + x_{ij}^2) \cdot (p_{ij,k}^2 + q_{ij,k}^2)}{U_{i,k}} = 0, \forall ij \in \Omega, k \in \mathbf{T} \quad (4)$$

$$p_{ij,k}^2 + q_{ij,k}^2 \leq s_{ij}^2, \forall ij \in \Omega, k \in \mathbf{T} \quad (5)$$

$$V_{i,k} = V^{Ref}, \forall i \in \Phi^{Sub}, k \in \mathbf{T} \quad (6)$$

$$V_i^{Lower} \leq V_{i,k} \leq V_i^{Upper}, \forall i \in \Phi \setminus \Phi^{Sub}, k \in \mathbf{T} \quad (7)$$

where

$$\begin{cases} P_{i,k}^{PV/EV} = \pi_{i,k}^{PV/EV} \cdot S_i^{PV/EV}, \forall i \in \Phi, k \in \mathbf{T} \\ Q_{i,k}^{PV/EV} = \tan \lambda^{PV/EV} \cdot \pi_{i,k}^{PV/EV} \cdot S_i^{PV/EV}, \forall i \in \Phi, k \in \mathbf{T} \\ \pi_{i,k}^{PV/EV} \in [0,1], \forall i \in \Phi^{PV/EV}, k \in \mathbf{T} \\ \pi_{i,k}^{PV/EV} = 0, \forall i \in \Phi \setminus \Phi^{PV/EV}, k \in \mathbf{T} \end{cases} \quad (8)$$

$$Q_{i,k}^{OT} = \tan \lambda_{i,k}^{OT} \cdot P_{i,k}^{OT}, \forall i \in \Phi, k \in \mathbf{T} \quad (9)$$

In the above model, superscripts PV, EV and OT refer to PV generation, EVCSs and other loads, respectively. $P_{i,k}$ and $Q_{i,k}$ respectively denote the active and reactive power at bus $i$ in time slot $k$. $p_{ij,k}$ and $q_{ij,k}$ are respectively the active and reactive powers through branch $ij$ in time slot $k$. $r_{ij}$, $x_{ij}$ and $s_{ij}$ are the resistance, the reactance and the apparent power flow limitation of branch $ij$, respectively. $\pi_{i,k}$ is the proportion of actual outputs (PV generation) or loads (EVCS) to the maximum capacity at bus $i$ in time slot $k$. $\lambda$ is the power factor angle. $V_{i,k}$ is the voltage of bus $i$ in time slot $k$, and $V_i^{Upper}$ and $V_i^{Lower}$ are respectively the upper and lower bounds of the voltage of bus $i$. $V^{Ref}$ is the reference voltage for substation buses. $\Phi$ is the set of all the buses and $\Phi^{Sub}$ is the set of the substation buses ($\Phi^{PV}, \Phi^{EV}, \Phi^{Sub} \subseteq \Phi$). $\Omega_i$ is the set of branches connected to bus $i$ and $\Omega$ is the set of all the branches. $\mathbf{T}$ is the set of all the time slots concerned.

In constraints (2)-(7), (2) and (3) are active and reactive power balance constraints; (4) are distribution system power flow equations for radial networks [24]; (5) describe the branch capacity limitation; (6) and (7) are voltage constraints. To linearize the model, the nonconvex constraints (4) can be rewritten as (10) by neglecting the relatively small quadratic term of line losses and regarding $V_{i,k}^2$ as a new variable $U_{i,k}$, $\forall i \in \mathbf{\Phi}, k \in \mathbf{T}$ (same substitutions in (6) and (7) keep the constraints linear), and the circular constraints (5) can be linearized as (11) through polytopic approximation [25]. In (11), a circumscribed octagon is appropriately used to do the linearization.

$$U_{i,k} - U_{j,k} - 2 \cdot (r_{ij} \cdot p_{ij,k} + x_{ij} \cdot q_{ij,k}) = 0, \forall ij \in \mathbf{\Omega}, k \in \mathbf{T} \quad (10)$$

$$\begin{cases} |p_{ij,k}| \leq s_{ij}, |q_{ij,k}| \leq s_{ij} \\ |p_{ij,k} + q_{ij,k}| \leq \sqrt{2} \cdot s_{ij}, \forall ij \in \mathbf{\Omega}, k \in \mathbf{T} \\ |p_{ij,k} - q_{ij,k}| \leq \sqrt{2} \cdot s_{ij} \end{cases} \quad (11)$$

By replacing (4) and (5) with (10) and (11), an approximate linear model is then obtained.

*B. Description of Uncertainty*

In light of the model built hereinbefore, we introduce the uncertainties in this subsection. Random variables $\tilde{\xi}_{i,k}^{\text{PV}}$, $\tilde{\xi}_{i,k}^{\text{EV}}$ and $\tilde{\xi}_{i,k}^{\text{OT}}$ are used to describe the output or load uncertainties of PV generation, EVCSs and other loads. Then $P_{i,k}^{\text{PV/EV/OT}}$ should be replaced by $\tilde{P}_{i,k}^{\text{PV/EV/OT}}$ whose expressions are as shown in (12).

$$\tilde{P}_{i,k}^{\text{PV/EV/OT}} = \left(1 + \tilde{\xi}_{i,k}^{\text{PV/EV/OT}}\right) \cdot P_{i,k}^{\text{PV/EV/OT}}, \forall i \in \mathbf{\Phi}, k \in \mathbf{T} \quad (12)$$

Note that 1) the uncertainties of PV generation and EVCSs are in fact reflected in $\pi_{i,k}^{\text{PV}}$ and $\pi_{i,k}^{\text{EV}}$; 2) power factors are supposed to be constant as $\tilde{\xi}_{i,k}^{\text{PV}}$, $\tilde{\xi}_{i,k}^{\text{EV}}$ and $\tilde{\xi}_{i,k}^{\text{OT}}$ vary so that the uncertainties of $\tilde{P}_{i,k}^{\text{PV/EV/OT}}$ and $\tilde{Q}_{i,k}^{\text{PV/EV/OT}}$ are coincident.

Considering that the random variables usually cannot vary in infinite domains in practice, uncertainty intervals with upper and lower bounds are necessary for more accurate characterization of the random vector $\tilde{\boldsymbol{\xi}} = \{\tilde{\xi}_{i,k}^{\text{PV}}, \tilde{\xi}_{i,k}^{\text{EV}}, \tilde{\xi}_{i,k}^{\text{OT}}\}$. Suppose dimension of $\tilde{\boldsymbol{\xi}}$ is $Z$, and $\tilde{\xi}_z \in [\xi_z^{\text{Lower}}, \xi_z^{\text{Upper}}]$, $\forall z = 1, \cdots, Z$. Then, the support information of $\tilde{\boldsymbol{\xi}}$ can be expressed as (13).

$$\boldsymbol{\Xi} = \left\{ \tilde{\boldsymbol{\xi}} \in \mathbb{R}^Z \middle| [\tilde{\boldsymbol{\xi}}^{\text{T}}, 1] \mathbf{W}_z [\tilde{\boldsymbol{\xi}}^{\text{T}}, 1]^{\text{T}} \leq 0, \forall z = 1, \cdots, Z \right\} \quad (13)$$

where

$$\mathbf{W}_z = \begin{bmatrix} & 1 & & -(\xi_z^{\text{Lower}} + \xi_z^{\text{Upper}})/2 \\ & & & \\ & -(\xi_z^{\text{Lower}} + \xi_z^{\text{Upper}})/2 & & \xi_z^{\text{Lower}} \cdot \xi_z^{\text{Upper}} \end{bmatrix} \begin{matrix} z \\ \\ Z+1 \end{matrix}$$

$z \qquad Z+1$

$\forall z = 1, \cdots, Z \quad (14)$

Further, we consider the uncertainty based on the first and second order moment information, which can be extracted from the historical data. Let $\boldsymbol{\mu} \in \mathbb{R}^Z$ be the mean vector and $\boldsymbol{\Sigma} \in \mathbb{S}^Z$ ($\mathbb{S}^Z$ denotes the space of symmetric matrices of dimension $Z$) be the covariance matrix of $\tilde{\boldsymbol{\xi}}$. Then, the set of probability distribution supported on $\boldsymbol{\Xi}$ with the same moment information can be expressed as (15).

$$\boldsymbol{\Theta}_{\boldsymbol{\Xi}} = \left\{ F \middle| \mathbb{E}_F(1) = 1, \mathbb{E}_F(\tilde{\boldsymbol{\xi}}) = \boldsymbol{\mu}, \mathbb{E}_F(\tilde{\boldsymbol{\xi}} \cdot \tilde{\boldsymbol{\xi}}^{\text{T}}) = \boldsymbol{\Sigma} + \boldsymbol{\mu} \cdot \boldsymbol{\mu}^{\text{T}} \right\} \quad (15)$$

where $\mathbb{E}_F(\cdot)$ is the expectation taken with respect to the $\tilde{\boldsymbol{\xi}}$, given that it follows the probability distribution $F$. In (15), $\mathbb{E}_F(1) = 1$ ensures the sum of the probabilities is 1, and $\mathbb{E}_F(\tilde{\boldsymbol{\xi}}) = \boldsymbol{\mu}$ and $\mathbb{E}_F(\tilde{\boldsymbol{\xi}} \cdot \tilde{\boldsymbol{\xi}}^{\text{T}}) = \boldsymbol{\Sigma} + \boldsymbol{\mu} \cdot \boldsymbol{\mu}^{\text{T}}$ enforce consistency with the given first and second order moments, respectively. The set $\boldsymbol{\Theta}_{\boldsymbol{\Xi}}$ makes full exploitation of the easily accessible information, i.e., the moments and the support set, to comprehensively describe the uncertainty of the uncertain parameters.

*C. Formulation of Distributionally Roubust Joint Constrained Programming Model*

We rewrite the linearized model built in Subsection II-A in form of (16) and (17).

$$\min \boldsymbol{c}^{\text{T}} \cdot \boldsymbol{x} \quad (16)$$

subject to

$$\boldsymbol{a}_m(\tilde{\boldsymbol{\xi}})^{\text{T}} \cdot \boldsymbol{x} \leq b_m(\tilde{\boldsymbol{\xi}}), \forall m = 1, \cdots, M \quad (17)$$

The objective function (16) is the vector form of (1), where $\boldsymbol{c} = \{\alpha_i^{\text{PV}}, \alpha_i^{\text{EV}}\}$ and $\boldsymbol{x} = \{S_i^{\text{PV}}, S_i^{\text{EV}}\}$ both with dimension $N$. In (17), $M$ is the number of constraints, and $\boldsymbol{a}_m(\tilde{\boldsymbol{\xi}}) \in \mathbb{R}^N$ and $b_m(\tilde{\boldsymbol{\xi}}) \in \mathbb{R}$ are uncertain constraint coefficients and can be expressed as below.

$$\boldsymbol{a}_m(\tilde{\boldsymbol{\xi}}) = \boldsymbol{a}_m^0 + \sum_{z=1}^{Z} \tilde{\xi}_z \cdot \boldsymbol{a}_m^z, \forall m = 1, \cdots, M, z = 1, \cdots, Z \quad (18)$$

$$b_m(\tilde{\boldsymbol{\xi}}) = b_m^0 + \sum_{z=1}^{Z} \tilde{\xi}_z \cdot b_m^z, \forall m = 1, \cdots, M, z = 1, \cdots, Z \quad (19)$$

For ease of notation in the following, we utilize auxiliary functions $y_m^j : \mathbb{R}^N \to \mathbb{R}$, as shown in (20), to rewrite (17) as (21).

$$y_m^z(\boldsymbol{x}) = (\boldsymbol{a}_m^z)^{\text{T}} \cdot \boldsymbol{x} - b_m^z, \forall m = 1, \cdots, M, z = 1, \cdots, Z \quad (20)$$

$$\boldsymbol{y}_m(\boldsymbol{x})^{\text{T}} \cdot \tilde{\boldsymbol{\xi}} + y_m^0(\boldsymbol{x}) \leq 0, \forall m = 1, \cdots, M \quad (21)$$

where $\boldsymbol{y}_m(\boldsymbol{x}) = [y_i^1(\boldsymbol{x}), \cdots, y_i^k(\boldsymbol{x})]^{\text{T}}$. However, constraints (21) are over-conservative, and to fix it, (21) are converted into a DR joint chance constraint (22).

$$\inf_{F \in \boldsymbol{\Theta}_{\boldsymbol{\Xi}}} \mathbb{P}_F \left( \boldsymbol{y}_m(\boldsymbol{x})^{\text{T}} \cdot \tilde{\boldsymbol{\xi}} + y_m^0(\boldsymbol{x}) \leq 0, \forall m = 1, \cdots, M \right) \geq 1 - \varepsilon \quad (22)$$

where $\mathbb{P}_F(\cdot)$ is the probability taken with respect to the $\tilde{\boldsymbol{\xi}}$, given that it follows the probability distribution $F$, and $1-\varepsilon$ is the confidence level ($\varepsilon$ reflects the risk tolerance) which can be used to control the risk. By replacing (17) with (22), a DR joint chance constrained programming model for risk-averse joint capacity evaluation of PV generation and EVCSs is realized.

## III. SOLUTION TO MODEL OF RISK-AVERSE JOINT CAPACITY EVALUATION

### A. Transformation into WC-CVaR Constraint

According to [26], constraint (22) can be further reformulated as a DR individual chance constraint, expressed in (23).

$$\inf_{F \in \Theta_{\boldsymbol{\xi}}} \mathbb{P}_F \left( \max_{m=1,\cdots,M} \left\{ \gamma_m \cdot \left( \boldsymbol{y}_m(\boldsymbol{x})^{\mathrm{T}} \cdot \tilde{\boldsymbol{\xi}} + y_m^0(\boldsymbol{x}) \right) \right\} \leq 0 \right) \geq 1-\varepsilon \quad (23)$$

where $\gamma_m$ is positive scaling parameters and does not affect the feasible region of (23), but can be used to improve the approximation to be developed below. Constraint (23) can be transformed into the conservative form of a WC-CVaR constraint [26], as (24).

$$\sup_{F \in \Theta_{\boldsymbol{\xi}}} \text{CVaR}_\varepsilon \left( \max_{m=1,\cdots,M} \left\{ \gamma_m \cdot \left( \boldsymbol{y}_m(\boldsymbol{x})^{\mathrm{T}} \cdot \tilde{\boldsymbol{\xi}} + y_m^0(\boldsymbol{x}) \right) \right\} \right) \leq 0 \quad (24)$$

For any fixed $\boldsymbol{\gamma} = \{\gamma_m\}$, it has been proved that (24) can be conservatively reformulated in terms of linear matrix inequalities (LMIs) [27], as (25)-(28).

$$\beta + \frac{1}{\varepsilon} \text{Tr}(\mathbf{B} \cdot \mathbf{H}) \leq 0, \mathbf{B} = \begin{bmatrix} \boldsymbol{\Sigma} + \boldsymbol{\mu} \cdot \boldsymbol{\mu}^{\mathrm{T}} & \boldsymbol{\mu} \\ \boldsymbol{\mu}^{\mathrm{T}} & 1 \end{bmatrix} \quad (25)$$

$$\tau_m \geq 0, \forall m = 0, 1, \cdots, M \quad (26)$$

$$\mathbf{H} + \sum_{z=1}^{Z} \tau_{0,z} \mathbf{W}_z \succeq \mathbf{0} \quad (27)$$

$$\mathbf{H} + \sum_{z=1}^{Z} \tau_{m,z} \mathbf{W}_z - \begin{bmatrix} \mathbf{0} & \frac{1}{2} \cdot \gamma_m \cdot \boldsymbol{y}_m(\boldsymbol{x}) \\ \frac{1}{2} \cdot \gamma_m \cdot \boldsymbol{y}_m(\boldsymbol{x})^{\mathrm{T}} & \gamma_m \cdot y_m^0(\boldsymbol{x}) - \beta \end{bmatrix} \succeq \mathbf{0},$$

$$\forall m = 1, \cdots, M \quad (28)$$

where $\mathbf{H} \in \mathbb{S}^N$, $\beta \in \mathbb{R}$, $\boldsymbol{\tau}_m \in \mathbb{R}^Z$, $\text{Tr}(\cdot)$ denotes the trace of the matrix, and $\mathbf{X} \succeq \mathbf{0}$ implies that matrix $\mathbf{X}$ is positive semidefinite. Actually, excluding decision variables $\boldsymbol{x}$ and auxiliary variables $\beta$ and $\mathbf{H}$, scaling parameters $\boldsymbol{\gamma}$ are also variables of the model composed of objective (16) and constraints (25)-(28). As a consequence, the model includes BMI constraints, which make the problem hard to solve. Leveraging the property that the model is equivalent to a tractable SDP for any fixed $\boldsymbol{\gamma}$, a method is presented as the solution in the next subsection.

### B. Iterative Algorithm Based on SDP

The algorithm introduced here is inspired by [26] and [27]. We first construct an auxiliary optimization problem:

$$\min \beta + \frac{1}{\varepsilon} \text{Tr}(\mathbf{B} \cdot \mathbf{H}) \quad (29)$$

subject to

$$(26)\text{-}(28) \quad (30)$$

For notation convenience, we call the model with objective (16) and constraints (25)-(28) *Model A* and the above model *Model B*. Both *Model A* and *Model B* are SDPs if $\boldsymbol{x}$ or $\boldsymbol{\gamma}$ is fixed. Based on this property, the problem can be solved by following iterative algorithm based on SDP:

*1)* For a given $\bar{\boldsymbol{\gamma}}$, an optimal $\hat{\boldsymbol{x}}$ can be obtain by solving *Model A*.

*2)* Let $\hat{\boldsymbol{x}}$ be frozen and solve *Model B* to obtain the optimal scaling parameters $\hat{\boldsymbol{\gamma}}$.

*3)* Let $\bar{\boldsymbol{\gamma}} = \hat{\boldsymbol{\gamma}}$ and repeat the above *1)* and *2)* until the decrement of objective (1) is less than a predetermined tolerance.

Due to the feasibility of $\hat{\boldsymbol{x}}$ and $\bar{\boldsymbol{\gamma}}$ in *Model A*, $\hat{\boldsymbol{x}}$ and $\hat{\boldsymbol{\gamma}}$ must be a feasible solution to *Model A*. As a result, a series of monotonically decreasing values of objective (1) can be obtained by repeating the above procedure.

## IV. CASE STUDIES

In this section, the proposed risk-averse joint capacity evaluation method is tested numerically on the IEEE 33-bus radial distribution system and the corresponding simulation results are presented and analyzed.

### A. Case Description

In this case, we let the weighting coefficients for capacities of PV generation and EVCSs be same. The test distribution system is shown in Fig. 1 with candidate buses of PV generation and EVCS marked, where PV generation candidate buses tend to be at the end of the radial networks for voltage support. The substation voltage magnitude (bus 1) is set as 10.5 kV and the permissible tolerance of bus voltages is [-0.5 kV, +0.3 kV]. The detailed parameters, including branch impedances, loads of buses and branch apparent power flow limits, of the test distribution system are listed in [28]. Besides, the load and PV output profiles are taken into account discretely by regarding the powers in each hour as a constant. The profiles of per-unit values of PV outputs, EV charging loads and other loads are depicted in Fig. 2. Given a bus, the basic values for the profiles are determined according to the corresponding load information (see [28]), PV generation capacity and EVCSs capacity, respectively. The power factors of PV outputs and EVCSs are set as 0.95 and 0.97, respectively. Regarding the uncertainties, for simplicity, we only consider the uncertainties for different types of loads/outputs, i.e., PV outputs, EV charging loads and other loads, while neglect the difference of uncertainties from bus to bus. Let the mean values, the variances and the lower and upper bounds of $\tilde{\xi}^{\text{PV/EV/OT}}$ be respectively $\mu^{\text{PV/EV/OT}}$, $d^{\text{PV/EV/OT}}$, $\xi^{\text{Lower,PV/EV/OT}}$ and $\xi^{\text{Upper,PV/EV/OT}}$. Since $\tilde{\xi}^{\text{PV}}$, $\tilde{\xi}^{\text{EV}}$ and $\tilde{\xi}^{\text{OT}}$ are independent of each other, $\boldsymbol{\Sigma}$ is a diagonal matrix with diagonal elements $d^{\text{PV}}$, $d^{\text{EV}}$ and $d^{\text{OT}}$.

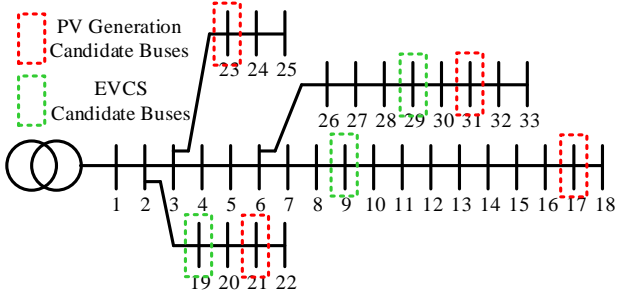

Fig. 1. The schematic diagram of IEEE 33-bus radial distribution system with candidate buses marked.

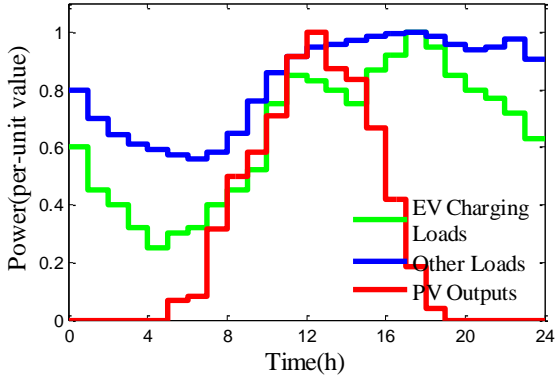

Fig. 2. The profiles of PV outputs, EV charging loads and other loads.

TABLE I.  PV GENERATION CAPACITIES AT CANDIDATE BUSES (MW)

|  | Bus 17 | Bus 21 | Bus 23 | Bus 31 |
|---|---|---|---|---|
| Capacity | 1.13 | 1.66 | 2.79 | 2.36 |

TABLE II.  EVCS CAPACITIES AT CANDIDATE BUSES (MW)

|  | Bus 9 | Bus 19 | Bus 29 |
|---|---|---|---|
| Capacity | 1.17 | 3.35 | 2.08 |

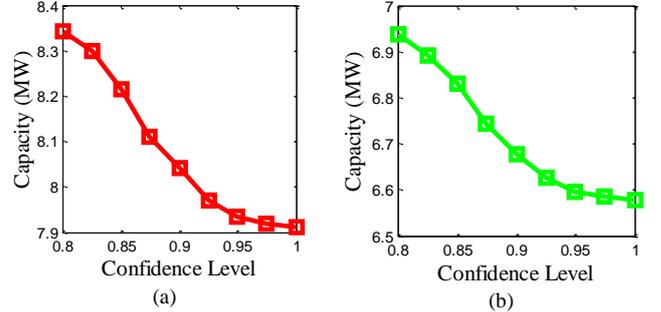

Fig. 3. The curves of maximum total capacities of (a) PV generation and (b) EVCSs versus confidence level.

## B. Results and Analysis

Let $\mu^{PV/EV/OT}$ be 0, $d^{PV/EV/OT}$ be 0.01, [$\xi^{Lower,PV/EV/OT}$, $\xi^{Upper,PV/EV/OT}$] be [-0.25, +0.25], and $1-\varepsilon$ be 0.95. The numerical results of capacities of PV generations and EVCSs at candidate buses are listed in Tables I and II, respectively. Supposing $1-\varepsilon$ varies, the capacity curves of PV generations and EVCSs are depicted in Fig. 3. When $\mu^{PV}$ and $\mu^{EV}$ can be changed, the surfaces of capacities are presented in Fig. 4. And the surfaces of capacities as $d^{PV}$ and $d^{EV}$ vary are presented in Fig. 5.

According to Tables I and II, it can be seen that a larger capacity of PV generation (EVCSs) of a bus can be achieved if there are EVCSs (PV generators) at a nearby bus and the connecting branch capacities are acceptable. Because adjacent locations of PV generators and EVCSs can mitigate the voltage deviation caused by long distance transmission and promote the PV generation absorption and EV integration.

Fig. 3 implies that the capacities of PV generation and EVCSs decline with the increasing of confidence level and the capacities converge to lower values as the confidence level tends to 1. Actually, when $1-\varepsilon=1$, all the constraints must always be met and the solutions are robust. From Fig. 4, it can be observed that 1) the total capacity increase slightly when $\mu^{PV}$ and $\mu^{EV}$ decrease; 2) the increasing of $\mu^{PV}$ ($\mu^{EV}$) has a positive effect on the capacity of EVCSs (PV generation), but reduces the capacity of PV generation (EVCSs). The reasons are that the higher $\mu^{PV}$ ($\mu^{EV}$) means the higher proportion of actual PV outputs (EV charging loads) to the capacity of PV generation (EVCSs), which reduces the capacity of PV generation (EVCSs), and the greater probability, i.e., the longer duration, of high PV outputs (EV charging loads), which calls for more loads (generation). However, the effect of the former is stronger than that of the latter so that the total capacity increases with the decreasing of $\mu^{PV}$ and $\mu^{EV}$. Similar explanations can be used to interpret the effects of second moments (variances $d^{PV}$ and $d^{EV}$) in Fig. 5.

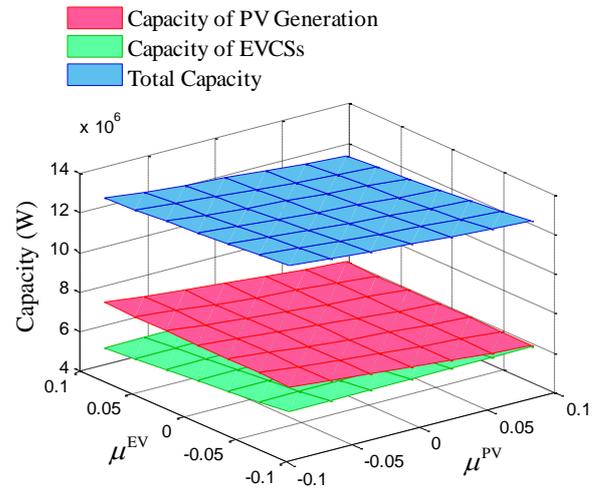

Fig. 4. The surfaces of capacities versus $\mu^{PV}$ and $\mu^{EV}$.

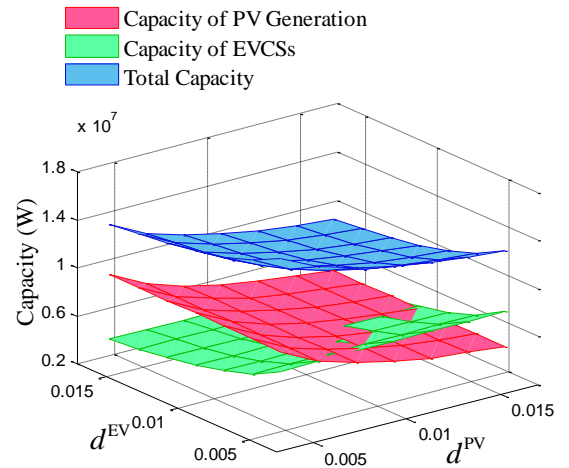

Fig. 5. The surfaces of capacities versus $d^{PV}$ and $d^{EV}$.

## V. CONSLUSIONS

In this work, we propose a risk-averse joint capacity evaluation method for PV generation and EVCSs, and a DR joint chance constrained programming model is accordingly formulated. The model leverages the limited and easily accessible information, i.e., moment and support information, of the uncertain parameters to determine the achievable joint capacity with an adjustable risk budget. Furthermore, the model can be transformed into a BMI constrained problem through WC-CVaR approximation and then be effectively solved by an iterative algorithm based on SDP. Case studies show that geographically close locations of PV generators and EVCSs can promote the insertion of them against the line congestions. And the analysis of the effects of random variable information on the capacities can provide guidelines for practical planning. In future work, the coordinated charging of EV will be considered in the capacity evaluation.

REFERENCES

[1] U.S. Environmental Protection Agency. (2015). *U.S. Greenhouse Gas Inventory Report: 1990-2014* [Online]. https://www.epa.gov/ghgemissions/us-greenhouse-gas-inventory-report-1990-2014, accessed Feb. 17, 2017.

[2] U.S. Energy Information Agency. (2015). *How much of the U.S. carbon dioxide emissions are associated with electricity generation?* [Online]. http://www.eia.gov/tools/faqs/faq.cfm?id=77&t=11, accessed Feb. 17, 2017.

[3] Y. Hua, M. Oliphant, E. J. Hu, "Development of renewable energy in Australia and China: A comparison of policies and status," *Renewable Energy*, vol. 85, pp. 1044-1051, Jan. 2016.

[4] S. Shinkhede, "Implementation of the Micro-Grid Concept and Balancing Massive Energy Production from Renewable Sources," *Int. Refereed J. Eng. and Sci.*, vol. 3, no. 1, pp. 76-84, Jan. 2014.

[5] International Energy Agency. (2015). *CO2 emissions from fuel combustion highlights 2015* [Online]. https://www.iea.org/publications/freepublications/publication/CO2EmissionsFromFuelCombustionHighlights2015.pdf, accessed Feb. 17, 2017.

[6] Z. Xu, Z. Hu, Y. Song, W Zhao and Y. Zhang, "Coordination of PEVs charging across multiple aggregators," *Appl. Energy*, vol. 136, pp. 582-589, Aug. 2014.

[7] L. Feng, S. Ge, and H. Liu, "Electric vehicle charging station planning based on weighted Voronoi diagram," in *Proc. Asia-Pac. Power Energy Eng. Conf.*, Shanghai, China, Mar. 2012, pp. 1-5.

[8] C. L. T. Borges, and D. M. Falcao, "Optimal distributed generation allocation for reliability, losses, and voltage improvement," *Int. J. Electr. Power Energy Syst.*, vol. 28, no. 6, pp. 413-420, July 2006.

[9] L. P. Fernández, T. G. S. Román, R. Cossent, C. M. Domingo, and P. Frías, "Assessment of the impact of plug-in electric vehicles on distribution networks." *IEEE Trans. Power Syst.*, vol. 26, no. 1, pp. 206-213, Feb. 2011.

[10] L. Dow, M. Marshall, and L. Xu, et al., "A novel approach for evaluating the impact of electric vehicles on the power distribution system," in *Proc. IEEE Power and Energy Society General Meeting*, Minneapolis, MN, USA, July 25-29, 2010.

[11] G. A. Putrus, P. Suwanapingkarl, and D. Johnston, *et al.*, "Impact of electric vehicles on power distribution networks," *IEEE Vehicle Power and Propulsion Conference*, Dearborn, MI, Sep. 7-10, 2009.

[12] D. Q. Hung and N. Mithulananthan, "Multiple distributed generators placement in primary distribution networks for loss reduction," *IEEE Trans. Ind. Electron.*, vol. 60, no. 4, pp. 1700-1708, Apr. 2013.

[13] V. R. Pandi, H. H. Zeineldin, and W. Xiao, "Determining optimal location and size of distributed generation resources considering harmonic and protection coordination limits," *IEEE Trans. Power Syst.*, vol. 28, no. 2, pp. 1245-1254, May. 2013.

[14] M. Raoofat, "Simultaneous allocation of DGs and remote controllable switches in distribution networks considering multilevel load model," *Int. J. Electr. Power Energy Syst.*, vol. 33, no. 8, pp. 1429-1436, Oct. 2011.

[15] A. M. El-Zonkoly, "Optimal placement of multi-distributed generation units including different load models using particle swarm optimisation," *IET Gener., Transm., Distrib.*, vol. 5, no. 7, pp. 760-771, Jul. 2011.

[16] F. Rotaru, G. Chicco, G. Grigoras, and G. Cartina, "Two-stage distributed generation optimal sizing with clustering-based node selection," *Int. J. Electr. Power Energy Syst.*, vol. 40, no. 1, pp. 120-129, Sep. 2012.

[17] A. H. Hajimiragha, C. A. Canizares, M. W. Fowler, S. Moazeni, and A. Elkamel, "A robust optimization approach for planning the transition to plug-in hybrid electric vehicles," *IEEE Trans. Power Syst.*, vol. 26, no. 4, pp. 2264-2274, Nov. 2011.

[18] G. Wang, Z. Xu, F. Wen, and K. Wong, "Traffic-constrained multiobjective planning of electric-vehicle charging stations," *IEEE Trans. Power Del.*, vol. 28, no. 4, pp. 2363-2372, Oct. 2013.

[19] M. F. Shaaban and E. F. El-Saadany, "Accommodating high penetrations of PEVs and renewable DG considering uncertainties in distribution systems," *IEEE Trans. Power Syst.*, vol. 29, no. 1, pp. 259-270, 2014.

[20] M. H. Moradi, M. Abedini, S. R. Tousi, and S. M. Hosseinian, "Optimal siting and sizing of renewable energy sources and charging stations simultaneously based on Differential Evolution algorithm," *Int. J. Elec. Power.*, vol. 73, pp. 1015-1024, 2015.

[21] D. Quoc, Z. Yang, and H. Trinh, "Determining the size of PHEV charging stations powered by commercial grid-integrated PV systems considering reactive power support," *Appl. Energy*, vol. 183, pp. 160-169, 2016.

[22] M. H. Amini, M. P. Moghaddam, and O. Karabasoglu, "Simultaneous allocation of electric vehicles' parking lots and distributed renewable resources in power distribution network," *Sustainable Cities and Society*, vol. 28, pp. 332–342, 2017.

[23] K. Zou, A. P. Agalgaonkar, K. M. Muttaqi and S. Perera, "Distribution system planning with incorporating DG reactive capability and system uncertainties," *IEEE Trans. Sustainable Energy*, vol. 3, no. 1, pp. 112-123, Jan. 2012.

[24] M. Baran and F. F. Wu, "Optimal sizing of capacitors placed on a radial distribution system," *IEEE Trans. Power Delivery*, vol. 4, no. 1, pp. 735-743, Jan 1989.

[25] C. N. Jones and M. Morari, "Polytopic approximation of explicit model predictive controllers," *IEEE Trans. Autom. Control*, vol. 55, no. 11, pp. 2542–2553, Nov. 2010.

[26] W. Chen, M. Sim, J. Sun, C-P. Teo, "From CVaR to uncertainty set: implications in joint chance-constrained optimization," *Oper. Res.*, vol. 58, no. 2, pp. 470-485, Mar. 2010.

[27] S. Zymler, D. Kuhn, and B. Ruste, "Distributionally robust joint chance constraints with second-order moment information," *Mathematical Programming*, vol. 137, no. 1, pp. 167-198, Feb. 2013.

[28] (Feb. 17, 2017). *IEEE 33-Bus Radial Distribution System* [Online]. https://drive.google.com/file/d/0B5QzF1wNoDzjZnl6REJuLWF0TTQ/view, accessed Feb. 17, 2017.